\documentstyle[12pt]{amsart}

\newtheorem{theorem}{Theorem}[section]
\newtheorem{lemma}[theorem]{Lemma}

\newtheorem{proposition}[theorem]{Proposition}

\theoremstyle{definition}

\newtheorem{example}[theorem]{Example}

\newtheorem{remark}[theorem]{Remark}

\newcommand{\cp}{{\Bbb C}_p}
\newcommand{\op}{{\mathcal O}_p}
\newcommand{\p}{{\mathcal P}}

\font\bb=msbm10 at 12pt

\newcommand{\bbQ}{\mbox{\bb Q}}

\begin{document}
\title{Morphic heights and periodic points}
\author{M. Einsiedler}
\address{(M.E.) Mathematical Institute, University of Vienna,
Strudlhofgasse 4, A-1090 Wien, Austria.}
\email{manfred@@mat.univie.ac.at}
\author{G. Everest}
\author{T. Ward}
\address{(G.E. \& T.W.) School of Mathematics, University of East Anglia,
Norwich NR4 7TJ, UK.}
\email{g.everest@@uea.ac.uk}
\email{t.ward@@uea.ac.uk}
\dedicatory{\today}
\thanks{The first author acknowledges the support of
EPSRC postdoctoral award GR/M49588, the
second thanks Jonathan Lubin and
Joe Silverman for the AMS Sectional meeting on Arithmetic
Dynamics at Providence, RI, 1999}
\begin{abstract}
An approach to the calculation of local
canonical morphic heights is described,
motivated by the analogy between the
classical height in Diophantine geometry and
entropy in algebraic dynamics.
We consider cases where the local morphic height is
expressed as an
integral average of the logarithmic distance to the
closure of the periodic points of the
underlying morphism.
The results may be thought
of as a kind of
morphic Jensen formula.
\end{abstract}
\maketitle
\section{Introduction}

Let $\phi:{\Bbb P}^1({\overline{\Bbb Q}})\to
{\Bbb P}^1({\overline{\Bbb Q}})$ be a morphism of degree $d$,
defined over the rationals.
Call, Goldstine and Silverman
(see \cite{call-goldstine-1997},\cite{call-silverman-1993})
have associated to $\phi$ a canonical
global {\sl morphic height}
$\hat{\lambda}_{\phi}$ on $\overline{\Bbb Q}$, with the
properties that
\begin{enumerate}
\item $\hat{\lambda}_{\phi}(\phi(q)) = d\hat{\lambda}_{\phi}(q)$
for any $q\in{\Bbb P}^1(\overline{\Bbb Q})$;
\item $q$ is pre-periodic if and only if $\hat{\lambda}_{\phi}(q) = 0$.
\end{enumerate}
A point $q$ is called {\sl pre-periodic} under
$\phi$ if the
orbit $\{\phi^{(n)}(q)\}$ is finite (write
$f^{(n)}$ for the $n$th iterate
of a map $f$).
The global height decomposes into a sum of local canonical
morphic heights $\lambda_{\phi,v}$:
\begin{equation*}
\hat{\lambda}_{\phi}(q)=\sum_{v}n_v\lambda_{\phi,v}(q).
\end{equation*}
Here $v$ runs over all the valuations (both finite
and infinite) of the
number field generated by $q$ and the $n_v$ denote the usual
normalising constants.
In the special case that $\phi$ takes the form
$\phi[x,y]=[y^df(x/y),y^d]$ for a polynomial
$f$ of degree $d$, Call and Goldstine
\cite{call-goldstine-1997} prove that
\begin{equation}\label{localheightisnicelimit}
\lambda_{\phi,v}(q)=\lim_{n\rightarrow \infty}
\frac{1}{d^n}\lambda_v(\phi^{(n)}(q)),
\end{equation}
where $\lambda_v$ is the local projective height
$\lambda_v(q)=\log ^+|q|_v$, and
that the local height $\lambda_{\phi,v}(q)$ vanishes if and only if
$|\phi^{(n)}(q)|_v$ is bounded for all $n$, and
finally, that $q$ is pre-periodic if and only if
\begin{equation}\label{preperiodifflocalheightvanish}
 \lambda_{\phi,v}(q)=0
\text{ for all $v$}.
\end{equation}

\begin{example}\label{paperbackwriter}
\begin{enumerate}
\item Let $f(z)=z^d$ with $d>1$. Here the canonical morphic
heights and the projective heights agree. Jensen's formula
(\cite[Theorem 15.18]{rudin-real-and-complex-analysis}) gives
$$
\int_{{\Bbb S}^1}\log |y-q|\mbox{d}m(y) = \log ^+ |q|,
$$
where $m$ is the Haar measure on the circle ${\Bbb S}^1$.
The circle is also the Julia set for this
morphism on ${\Bbb C}$, and it is the closure of the set of
non-zero periodic points,
which are all roots of unity. In the $p$-adic case, the Julia
set is empty but it is still true that the local height is the
Shnirel'man integral of the logarithmic distance
from the closure of the set of
periodic points. For all $v$, finite and infinite, the following holds:
$$
\log ^+ |q|_v =
\lim_{n\rightarrow \infty}\frac{1}{d^n}
\sum_{\zeta \neq q:\zeta^{d^n-1}=1}\log |\zeta -q|_v
$$
For finite $v$ this may be seen,
for instance, using Section \ref{reallythelastthing}.
Alternatively, it follows from the Diophantine estimate
$$
\vert q^{d^n}-q\vert_v>C(q)\vert n\vert_v
$$
provided the left hand side is non-zero
(cf. Remark \ref{omega}).
For $v\vert\infty$ a result from transcendence theory
(in this case, Baker's theorem) is needed.
\item Suppose $a,b\in \bbQ$ with $4a^3+27b^2\neq 0$ and let
$$
f(z)=\frac{z^4-2az^2-8bz+a^2}{z^3+az+b}.
$$
Then $f$ gives rise to a morphism of degree 4 which describes
the duplication map on an elliptic curve. The global and local
morphic heights coincide with the usual notions of height on the curve.
For the infinite valuations, the local height is again the
integral average
of the logarithmic distance to points on the Julia set, which
is
the closure of the set of periodic points.
At the
singular reduction primes it is still true that the local height
is the integral average of the logarithmic distance
to the periodic points. Both of these assertions are proved
in
\cite{einsiedler-everest-ward-2000} where this morphism was used
to construct a dynamical system which interprets these heights
in dynamical
terms.
The proofs require elliptic transcendence theory to show that
a rational point cannot approximate a periodic point too
closely (cf. Proposition \ref{lastthing}).
This is part of a much broader analogy
between heights in Diophantine
geometry and entropy in
algebraic dynamics (see
\cite{dambros-everest-miles-ward},
\cite{einsiedler-everest-ward-2000},
\cite{everest-ward-1998},
\cite{lind-ward-1988}).
\end{enumerate}
\end{example}

In this paper, our purpose is to describe
a family of examples
where the local canonical morphic height can be
expressed as a limiting integral over periodic points
of the underlying morphism.
The finite and infinite cases require different
approaches. In both, we consider
the special class of morphisms corresponding to affine polynomial
maps and in the former case, we assume good reduction in the sense of
Morton and Silverman
\cite{morton-silverman-1995}.

There are two directions in which this work
can be made more sophisticated that are not pursued here.
The first is to give a more formal interpretation of
the limiting process using Shnirel'man integrals
(see \cite{shnirelman}, or
\cite{koblitz} for a modern treatment); the second
is to extend the arguments to other morphisms.

\section{Complex case}

Assume that $\phi:{\Bbb P}^1({\Bbb C})\to{\Bbb P}^1({\Bbb C})$ is a morphism
of degree $d$, with
$\phi[x,y]=[y^df(x/y),y^d]$ for a polynomial
$f$ of degree $d$. For basic definitions
of complex dynamics, consult \cite{beardon-1991}.
The following theorem expresses the
local morphic height as an integral over the Julia
set $J(f)$  of the polynomial $f$,
and standard results from
complex dynamics show that this is in turn a limiting
integral over periodic points.

\begin{theorem}\label{doughnut}
If $f(z)=az^d+\dots+a_0$ is a polynomial,
then for any $q\in{\Bbb C}$,
\begin{equation}\label{juliaintegral}
\lambda_{\phi,\infty}(q)=\frac{1}{d-1}
\log\vert a\vert+\int_{J(f)}\log\vert x-q\vert\mbox{\textup d}m(x),
\end{equation}
where $m$ is the maximal invariant measure for $f$ on $J(f)$.
\end{theorem}

\begin{pf} Assume first that $q$ is in the domain of attraction
of $\infty$ for $f$.
The zeros of the
polynomial $f_n(x)=f^{(n)}(x)-x$ are precisely the solutions of the
equation $f^{(n)}(x)=x$.
Note that $d_n=\deg(f_n)=d^n$, where $d=\deg(f)$.
Since $|f^{(n)}(q)|\to\infty$,
$\frac{1}{d_n}\log\vert f_n(q)\vert$ is approximately
$\frac{1}{d^n}\log\vert f^{(n)}(q)\vert$,
which converges to $\lambda_{\phi,\infty}(q)$.
Since $q$ lies in the open Fatou set,
$\log|x-q|$ is continuous on $J(f)$. Now
\begin{equation}
\label{iwasarrested}
\frac{1}{d_n}\log|f_n(q)|=\frac{1}{d_n}
\sum_{f^{(n)}(x)=x}\log|x-q|+\frac{1}{d_n}\log\vert B_n\vert,
\end{equation}
where the sum is over the $n$th `division points'
and
\begin{equation*}
B_n=a^{1+d+d^2+\dots +d^{(n-1)}}
\end{equation*}
is the
leading coefficient of $f^{(n)}(x)$.
Thus
\begin{equation}
\label{holdingpapers}
\frac{1}{d_n}\log\vert B_n\vert=
\frac{1}{d^n}\left(\frac{d^n-1}{d-1}\right)
\log\vert a\vert\to\frac{1}{d-1}\log\vert a\vert.
\end{equation}
Now it is known that
$$\frac{1}{d_n}\sum_{f^{(n)}(x)=x}\log|x-q|\to\int_{J(f)}
\log|x-q|\text{d}m(x),$$
where $m$ is the maximal invariant measure for $f$ restricted to the Julia
set (see \cite{lyubich-1982};
\cite{MR85m:58110b}).

In the following it is convenient to assume that a=1, we can ensure this by
conjugating by a linear map.
Now $\vert x-f(q)\vert=\prod_{f(t)=x}\vert t-q\vert$, so
\begin{eqnarray}
\int_{J(f)}\log\vert x-f(q)\vert\mbox{d}m(x)&=&
\int_{J(f)}\sum_{f(t)=x}\log\vert t-q\vert\mbox{d}m(x)\nonumber\\
&=&d\int_{J(f)}\log\vert x-q\vert\mbox{d}m(x)\label{bagel}
\end{eqnarray}
(the last equality follows from \cite{lyubich-1982}
or \cite[Theorem (d)]{MR85m:58110b}).

Let now $q\notin J(f)$ have bounded orbit. Since $J(f)$ is closed, we can
find $\epsilon>0$ such that $B_{\epsilon}(q)\cap J(f)=\emptyset$. If
$\vert f^n(q)-q\vert>\epsilon/2$ for almost all $n$, then
\[
\frac{1}{d_n}\log\vert f^n(q)-q\vert\rightarrow 0.
\]
We can argue as in the first case to get
\[
 \int_{J(f)}\log\vert x-q\vert=0=\lambda_{\phi,\infty}(q).
\]
Assume now $\vert f^{n_j}(q)-q\vert\leq\epsilon/2$ for some sequence
$n_j\rightarrow\infty$. Then $\vert f^{n_j}(q)-x\vert>\epsilon/2$ for
$x\in J(f)$.
However, since $J(f)$ and $f^{n_j}(q)$ are bounded,
we have also an upper bound
\[
 \log\frac{\epsilon}{2}\leq\int_{J(f)}\log\vert x-f^{n_j}(q)\vert\leq M.
\]
Together with Equation (\ref{bagel}) we get
\[
 1/d^{n_j}\log\frac{\epsilon}{2}\leq\int_{J(f)}\log\vert x-q\vert
    \leq 1/d^{n_j}M
\]
which concludes the proof in this case.

It remains to show that the formula holds for $q\in J(f)$.
Since $J(f)$ has no interior, there is a sequence
$q_n\to q$ with $q_n\notin J(f)$.
Then $\log\vert x-q_n\vert\to
\log\vert x-q\vert$ for all $x\in J(f)\backslash\{q\}$.
Since $J(f)$ is bounded, $\log\vert x-q_n\vert$
and $\log\vert x-q\vert$ are uniformly bounded above
by $M$ say for $x\in J(f)\backslash\{q\}$.
So by Fatou's lemma
\begin{equation}
\label{bangbangmaxwell}
0=\lim_{n\to\infty}\int_{J(f)}\log\vert x-q_n\vert\mbox{d}m(x)
\le
\int_{J(f)}\log\vert x-q\vert\mbox{d}m(x)\le M.
\end{equation}
This shows that $x\mapsto\log\vert x-q\vert$ is
in $L^1(m)$.
If
$$
\int_{J(f)}\log\vert x-q\vert\mbox{d}m(x)>0,
$$
then
Equation (\ref{bagel}) contradicts (\ref{bangbangmaxwell}).
\end{pf}

\section{The $p$-adic case}
\label{reallythelastthing}

Let $\cp$ denote the usual completion of the algebraic
closure of the $p$-adic numbers ${\Bbb Q}_p$, and use
$\vert\cdot\vert$ to denote the extension of the $p$-adic
norm to $\cp$. Write $\op$ for the ring of integral
elements in $\cp$, and $\p$ for the maximal ideal
of $\op$.
In this section we
assume that $\phi:{\Bbb P}^1(\cp)\to{\Bbb P}^1({\Bbb C}_p)$ is a morphism
of degree $d$ corresponding to an
affine polynomial
$f$ of degree $d$ with coefficients in $\op$ and
leading coefficient in $\op^{*}$.
Notice that these assumptions are, for polynomials, equivalent to the
assumption that the map $\phi$
has good reduction in the sense of
\cite{morton-silverman-1995}: $\phi$ induces a morphism of
schemes over $\mbox{Spec}(\op)$.
The Julia set is empty in this setting (see
\cite{benedetto}, \cite{morton-silverman-1995}),
so a direct analogue of (\ref{juliaintegral}) is
not possible.

The main result expresses the local canonical
morphic height as a limiting integral over periodic points
for the polynomial $f$.

\begin{theorem}\label{maxwell}
If $\phi$ has good reduction and is defined by a polynomial
$f$ of degree $d$, then
$$
\lambda_{\phi,p}(q)=\log^{+}\vert q\vert=
\lim_{n\to\infty}\frac{1}{d^n}\sum_{\xi\neq q:
f^{(n)}(\xi)=\xi}\log\vert\xi-q\vert
$$
where the sum is taken with multiplicities.
\end{theorem}

Notice first that for $q\in\cp\backslash\op$ this
is clear, so from now on we assume that
$q\in\op$.
Despite the simple resulting value of the height, the
convergence involved requires an argument.
The main issue is to produce lower bounds on
the size of $\vert\xi-\zeta\vert$ for distinct periodic
points $\xi$ and $\zeta$.

The {\sl least period} of
a periodic point $\xi$ is the cardinality
of the orbit of $\xi$. The points of period $n$
are the solutions to the polynomial equation
\begin{equation}\label{onhishead}
f^{(n)}(x)-x=0,
\end{equation}
and are therefore all elements of $\op$.
Following \cite{morton-silverman-1995}, for
a periodic point $\xi$ define
$a_n(\xi)$ to be the multiplicity of $\xi$ in
(\ref{onhishead}), with the obvious convention that
$\xi$ has multiplicity zero in an equation that it
does not satisfy.
Notice that
$a_n(\xi)\neq 0$ if and only if $n$ is a multiple
of the least period of $\xi$.
Define $a_n^{*}(\xi)$
by
$$
a_n^{*}(\xi)=\sum_{d\vert n}\mu\left(\frac{n}{d}\right)
a_n(\xi),
$$
where $\mu$ is the M{\"o}bius function.
Increases in the multiplicity of the periodic
point $\xi$ along the sequence of
multiples of its least period
are recorded by $a_n^{*}(\xi)$.
The periodic point $\xi$ is an
{\sl essential $n$-periodic point} if
$a_n^{*}(\xi)>0$.

The following proposition is a special
case of \cite[Prop. 3.2]{morton-silverman-1995}.

\begin{proposition}\label{camedown}
Let $K$ be an algebraically closed field
of characteristic $p\ge0$,
and $f$ a polynomial over $K$ with degree $d\ge2$.
Fix a periodic point $\zeta\in K$ with
least period $m$, and let
$r$ denote the multiplicative order of
$(f^{(m)})^{\prime}(\zeta)$ in $K^{*}$
or $\infty$ if $(f^{(m)})^{\prime}$ is not a root of
unity.
Then for $n\ge 1$, $a_n^{*}(\zeta)\ge 1$ if and only if
one of the following conditions hold.
\begin{enumerate}
\item $n=m$;
\item $n=mr$;
\item $p>0$ and $n=p^emr$ for some $e\ge1$.
\end{enumerate}
\end{proposition}

When $K=\cp$, this proposition will also be applied to the polynomial
$\bar{f}$ induced by reduction mod $\p$.
Notice that the sum of the multiplicities of
the points of period $n$ under $f$ lying in one
residue class gives the multiplicity of the
image point as a point of period $n$ under $\bar{f}$.

For the proof of Theorem \ref{maxwell} the following
proposition is needed; this will be proved later.

\begin{proposition}\label{silverhammer}
Suppose $\xi$ is a periodic point with
least period $n$
for a polynomial $f$ of
good reduction. Then for any fixed $q \in\op$,
$\vert q-\xi\vert\to 1$ as $n\to\infty$, provided $q\neq \xi$.
\end{proposition}

\begin{pf}(of Theorem \ref{maxwell})
By Proposition \ref{camedown}
applied to the field $\cp$ with
characteristic zero,
if $\xi$ is
a periodic point with least period
$m$, then the multiplicity of $\xi$ viewed
as a periodic point of period $m,2m,3m,\dots$
is uniformly bounded.
Fix $q\in\op$ and $s\in(0,1)$.
Proposition \ref{silverhammer} says that
the number of periodic points
in the metric ball $D_s(q)$ is finite.
It follows that
$$
\liminf_{n\to\infty}\frac{1}{d^n}\sum_{\xi\neq q:
f^{(n)}(\xi)=\xi}\log\vert\xi-q\vert\ge\log s.
$$
On the other hand, each term in the sum
is non-positive, so letting $s\to 1$
proves the theorem.
\end{pf}

All that remains is to prove Proposition \ref{silverhammer},
for which we need some lemmas.

\begin{lemma}\label{wealllive}
Assume that $f(0)=0$,
and let $\zeta$ be a periodic point of $f$.
Then $\vert\zeta\vert=\vert f^{(n)}(\zeta)\vert$ for
all $n\ge1$.
\end{lemma}

\begin{pf}
The spherical metric
used in \cite{morton-silverman-1995} coincides with the usual metric
in $\op$, and $f$ has good reduction.
So by \cite[Prop. 5.2]{morton-silverman-1995},
$$
\vert f(x)-f(y)\vert\le\vert x-y\vert
$$
for $x,y\in\op$.
The lemma follows at once.
\end{pf}

\begin{lemma}\label{inayellowsubmarine}
Assume that $f(0)=0$ and $n>1$ is fixed.
\begin{enumerate}
\item If $\vert f^{\prime}(0)\vert<1$ then
$\prod_{\xi\neq0}\vert\xi\vert^{a_n^{*}(\xi)}=1$.
\item If $\vert f^{\prime}(0)-1\vert<p^{-1}$
then
$\prod_{\xi\neq0}\vert\xi\vert^{a_n^{*}(\xi)}=\left\{
\begin{array}{cl}
1/p&\mbox{if }n\mbox{ is a power of }p,\\
1&\mbox{if not.}
\end{array}
\right.$
\end{enumerate}
\end{lemma}

\begin{pf}
In case 1., $\bar{f}^{\prime}(0)=0$ in the algebraically
closed field
$\op/\p$, so Proposition \ref{camedown} may be applied with
$\zeta=0+\p$, $m=1$ and
$r=\infty$. It follows that for $n>1$ $a_n^{*}(0+\p)=0$,
so there cannot be an essential $n$-periodic point $\xi$
for $f$ with $\vert\xi\vert<1$.

In case 2., $m=1$ and $r=1$ for the point
$\zeta=0+\p$ in Proposition \ref{camedown}.
It follows that only values of $n$ of the form
$p^k$ are relevant.
Notice that
\begin{equation}\label{herecomesthesun}
\prod_{\xi\neq0}\vert\xi\vert^{a_{p^k}^{*}(\xi)}=
\Biggl\vert
\left(\frac{f^{(p^k)}(x)-x}{f^{(p^{k-1})}(x)-x}\right)_{x=0}
\Biggr\vert
\end{equation}
If $f^{\prime}(0)\neq1$, then
the right-hand side of
(\ref{herecomesthesun}) is given by
$$
\left\vert\frac{(f^{\prime}(0))^{p^k}-1}{(f^{\prime}(0))^{p^{k-1}}-1}
\right\vert=\frac{1}{p}
$$
by the binomial theorem.
If $f^{\prime}(0)=1$, write
$f(x)=x+x^eg(x)$
with $e>1$ and $g(0)\neq 0$, then a simple induction argument
shows that
$$
f^{(k)}(x)=x+kx^eg(x)+O(x^{2e-1}).
$$
It follows that (\ref{herecomesthesun}) is equal to $p$
again.

\end{pf}

\begin{pf}(of Proposition \ref{silverhammer})
Let $\zeta$ be any periodic point, with least period $\ell$.
The first step is to prove the proposition for $q=\zeta$.
Let $\xi$ have least period $n$ under $f$.
The multiplicity of $\zeta+\p$, which has least
period $m$ for some $m\vert\ell$ must increase
at $\ell$ (because $a_{\ell}^{*}(\zeta)>1$).
It follows by Proposition \ref{camedown} that
$\ell$ is equal to $m$, $mr$, or $mrp^e$ for some
$e\ge1$.
Assume first that $n$ is not of one of those
forms; then $\vert\xi-\zeta\vert=1$ because the
multiplicity of $\zeta+\p$ cannot increase
at $n$ in $\op/\p$.

In the remaining cases, we may assume for large
$n$ that $\ell\vert n$. Then $\xi$ is a periodic point
with least period $n/\ell$ under $f^{(\ell)}$.
Applying the conjugation $x\mapsto x-\zeta$ means that
$0$ is a fixed point of $g$, where $g(x)=f^{(\ell)}(x+\zeta)-
\zeta$.

If $\vert g^{\prime}(0)\vert<1$, then by Lemma \ref{inayellowsubmarine}
applied to $g$, $\vert\xi-\zeta\vert=1$.

If $\vert g^{\prime}(0)\vert=1$, let $t$ be the
order of $g^{\prime}(0)+\p$ in $\op/\p$.
Then
$$\vert(g^{\prime}(0))^{t}-1\vert<1.$$
There
exists a $c\ge1$
such that
$$\vert(g^{\prime}(0))^{tp^{c}}-1\vert<1/p.$$
As before, we may assume that $tp^c\vert n$, so
Lemma \ref{inayellowsubmarine}
may be applied to the map $h=g^{(tp^c)}$ to give
$$
\prod_{j=1,\dots,n/tp^c\ell}\vert h^{(j)}(\xi-\zeta)\vert\ge
1/p.
$$
It follows by Lemma \ref{wealllive} that
$\vert\xi-\zeta\vert\ge
p^{-tp^c\ell/n}$.
Since $t,c,\ell$ depend only on $\zeta$,
$\vert\xi-\zeta\vert\to1$
as $n\to\infty$. The ultrametric inequality in $\cp$ now gives the result
for any $q\in\op$.
\end{pf}

\begin{remark}\label{omega}
Notice that the discussion above also gives a
quantitative version of Proposition \ref{silverhammer}.
This Diophantine result
may be of independent interest.
If $f$ is a polynomial of good reduction, then
$$
\vert f^{(n)}(q)-q\vert> C(f,q)\vert n\vert
$$
for all $n\ge 1$, provided the left hand side is non-zero.
\end{remark}

\begin{example} To see the different cases that are
possible in Proposition \ref{silverhammer}, consider the
following examples.
\begin{enumerate}
\item Let $f(x)=g(x^p)+ph(x)$ be a monic polynomial
with coefficients in $\op$. Then
$\vert f^{\prime}(q)\vert<1$ for any $q\in\op$, so
in Lemma \ref{inayellowsubmarine} only the first
case is ever used. It follows that in Proposition \ref{silverhammer},
$\vert\zeta-\xi\vert=1$ for any distinct periodic points
$\zeta,\xi$.
\item Let $f(x)=x^2-(1+a)x$ for some small $a$.
Then $0$ is a fixed point, and
$\vert f^{\prime}(0)\vert=1$.
Now
$$
f^{(2)}(x)-x=x^4-(2a+2)x^3 + (a^2 + a)x^2 + (a^2 + 2a)x,
$$
so $(f^{(2)}(x)-x)/(f(x)-x)$ has constant term
$(a^2+2a)/(-a-2)=-a$.
Therefore there must be two non-zero points of
period $2$ that are close to the fixed point $0$.
\end{enumerate}
\end{example}

If the polynomial $f$ has coefficients
outside $\op$, then in contrast to Proposition \ref{silverhammer},
there may be sequences of periodic points
converging to a periodic point. For example,
$f(x)=x^2+\frac{1}{2}x$ on ${\Bbb C}_2$
has this property. Therefore, to recover Theorem \ref{maxwell}
in greater generality (for polynomials of
bad reduction or rational functions)
some kind of Diophantine approximation
results are needed.
In
Example \ref{paperbackwriter}.2
these tools are provided by
elliptic transcendence theory.

\begin{proposition}\label{lastthing}
 Let $\zeta$ be a periodic point with least
period $\ell$ under $f$.
 Assume that
$\vert (f^{(\ell)})^{\prime}(\zeta)\vert >1$. Then
there are periodic points $\xi\neq\zeta$
 arbitrarily close to $\zeta$.
\end{proposition}

\begin{pf}
 Define $g=f^{(\ell)}$ and $a= (f^{(\ell)})^{\prime}(\zeta)$.
 Without loss of generality we can assume that $\zeta=0$.
 Then
\[
 g(x)=ax+bx^e+O(x^{e+1})\mbox{ with }b\neq0
\]
and
\[
 g^{(2)}(x)=a^2x+(ab+a^eb)x^e+O(x^{e+1}).
\]
Define $b_2=(ab+a^eb)$, then $\vert b_2\vert=\vert b\vert\vert a\vert^e$.
 By induction one can see that
\[
 g^{(k)}(x)=a^kx+b_kx^e+O(x^{e+1})
\]
with $\vert b_k\vert=\vert b\vert \vert a\vert^{ke}$.

Therefore the Newton polygon of $g^{(k)}(x)-x$ starts with
a line with slope $s\leq -k+c$ for a fixed $c$ (depending on $b$
and $e$).
So there exists a periodic point $\xi$ with
$\vert \xi\vert= p^{s}\le p^{-k+c}$.
\end{pf}

\section{Tchebycheff polynomials}

\begin{example}\label{tchebysheffexample}
Consider
the Tchebycheff polynomial
of degree $d$,
$f(z)=T_d(z)=\cos(d\arccos(z))$.
The Julia set is the interval $J(f)=[-1,1]$.
The map $\phi:{\Bbb C}\to{\Bbb C}$ given by
$\phi(z)=\frac{1}{2}(z+z^{-1})$ is a semi-conjugacy from
$g:z\mapsto z^d$ onto $z\mapsto f(z)$, in other words,
$f(\phi(z))=\phi(z^d)$.
Write $\psi$ for the branch of the inverse of
$\phi$ defined on $\{z\in{\Bbb C}\mid\vert z\vert>1\}$.
The canonical morphic height at
the infinite place is (for $q\notin J(f)$)
\begin{eqnarray*}
\lambda_{\phi,\infty}(q)&=&\lim_{n\to\infty}\frac{1}{d^n}
\log^{+}\vert f^{(n)}(q)\vert\\
&=&\lim_{n\to\infty}\frac{1}{d^n}
\log^{+}\vert\phi g^{(n)}\psi(q)\vert\\
&=&\lim_{n\to\infty}\frac{1}{d^n}
\log^{+}\vert\frac{1}{2}\left(g^{(n)}\psi(q)+
\frac{1}{g^{(n)}\psi(q)}\right)\vert\\
&=&\lim_{n\to\infty}\max\left\{0,
\frac{1}{d^n}\log\vert g^{(n)}\psi(q)\vert\right\}\\
&=&\log^{+}\vert\psi(q)\vert.
\end{eqnarray*}
For $q\in J(f)$, the same formula holds since
${\lambda}_{\phi,\infty}(q)=0$ there by \cite{call-goldstine-1997}
and $\log^{+}\vert\psi(q)\vert=0$ there by
a direct calculation.

Now by Jensen's formula, for any $q\in{\Bbb C}$,
\begin{eqnarray*}
\log^{+}\vert\psi(q)\vert&=&
\log2+\int_{{\Bbb S}^1}\vert
\phi(y)-q\vert\mbox{d}y\\
&=&
\log2+\int_{J(f)}\log\vert t-q\vert\mbox{d}m(t)
\end{eqnarray*}
since $m$ is the image under $\phi$ of the
maximal measure (Lebesgue) on the circle.
That is,
\begin{equation}\label{staronpage5}
{\lambda}_{\phi,\infty}(q)=\log2+\int_{J(f)}\log\vert t-q\vert\mbox{d}m(t).
\end{equation}

The constant $\log 2$ in ${\lambda}_{\infty}(q)$
may be explained in accordance with
Theorem \ref{doughnut}. The
leading coefficient of $T_d$ is $2^{d-1}$,
so $\frac{1}{d-1}\log\vert a\vert$ in this case is
exactly $\log2$.
\end{example}

A similar approach can be adopted in the case
of polynomials with connected Julia sets.
There the local conjugacy near $\infty$
extends to the whole domain of attraction of $\infty$,
which is the complement of the filled Julia set.

\begin{example}\label{end}
As before, let $f(x)=T_d(x)=\cos(d\arccos(x))$ be the Tchebycheff polynomial
of degree $d$ and
let $\phi$ be the corresponding morphism.
We would like to use Theorem \ref{maxwell}, but
$f$ does not satisfy the assumptions since it is not monic.

Let $g(x)=2f(\frac{x}{2})$.
Notice that $f$ is defined uniquely by the property
$f(\frac{1}{2}(z+z^{-1}))=\frac{1}{2}(z^d+z^{-d})$.
It follows that $g$ is
characterized
by the property $g(z+z^{-1})=(z^d+z^{-d})$,
which shows that $g\in{\mathbb Z}[x]$ is a monic polynomial.
Let $\psi$ be the morphism defined by $g$,
then by Theorem \ref{maxwell} we have for $q\in{\mathbb C}_p$
\begin{equation}\label{notasong}
\lambda_{\psi,p}(q)=\log^{+}\vert q\vert=
\lim_{n\to\infty}\frac{1}{d^n}\sum_{\xi\neq q:
g^{(n)}(\xi)=\xi}\log\vert\xi-q\vert.
\end{equation}
Since $g(x)=2f(\frac{x}{2})$, we have that
$\lambda_{\phi,p}(q)=\lambda_{\psi,p}(2q)$
and on the right hand side of (\ref{notasong})
that $f^{(n)}(\xi)=\xi$ if and only if
$g^{(n)}(2\xi)=2\xi$. Therefore
\[
\lambda_{\phi,p}(q)=\log^{+}\vert 2q\vert=
\log\vert 2\vert+\lim_{n\to\infty}\frac{1}{d^n}\sum_{\xi\neq q:
f^{(n)}(\xi)=\xi}\log\vert\xi-q\vert,
\]
which is again analogous to Equation (\ref{staronpage5})
in Example \ref{tchebysheffexample}.
\end{example}

Example \ref{end} works because the Tchebycheff polynomial
can be conjugated to a polynomial of good reduction;
a similar approach can be adopted for any polynomial
that is conjugate to one of good reduction.

\bibliographystyle{amsplain}
%\bibliography{../bib/nt,../bib/etds,../bib/analysis,../bib/algebra}

\providecommand{\bysame}{\leavevmode\hbox to3em{\hrulefill}\thinspace}

\end{document}